\documentclass[a4paper, 12pt]{article}

\usepackage{amsfonts, amsmath, amsthm}
\usepackage{latexsym}

\theoremstyle{plain}
\newtheorem{thm}{Theorem}
\newtheorem{prp}{Proposition}
\newtheorem{lem}[prp]{Lemma} 

\theoremstyle{definition}

\theoremstyle{remark}
\newtheorem{rmk}{Remark}

\numberwithin{equation}{section}

\newcommand{\R}{\mathbb{R}}
\newcommand{\C}{\mathbb{C}}

\newcommand{\pa}{\partial}
\newcommand{\eps}{\varepsilon}

\DeclareMathOperator{\realpart}{\rm Re}
\DeclareMathOperator{\imagpart}{\rm Im}

\begin{document}
\title{\Large 
 Small data blow-up for a system of  nonlinear Schr\"odinger equations 
}  

\author{
         Tohru Ozawa\thanks{
             Department of Applied Physics, Waseda University. 
             3-4-1 Okubo, Shinjuku, Tokyo 169-8555, Japan. 
              (E-mail: {\tt txozawa@waseda.jp})}
         \and  
         Hideaki Sunagawa\thanks{
             Department of Mathematics, Graduate School of Science, 
             Osaka University. 
             1-1 Machikaneyama-cho, Toyonaka, Osaka 560-0043, Japan.
             (E-mail: {\tt sunagawa@math.sci.osaka-u.ac.jp})}
} 
 
\date{\today }   
\maketitle

\noindent{\bf Abstract:}\ We give examples of small data blow-up 
for a three-component system of quadratic nonlinear Schr\"odinger equations 
in one space dimension. Our construction of the blowing-up 
solution is based on the Hopf-Cole transformation, which allows us to reduce 
the problem to getting suitable growth estimates for a solution to 
the transformed system. Amplification in the reduced system is 
shown to have a close connection with the mass resonance.\\


\\


\section{Introduction}
We consider the initial value problem for a system of nonlinear 
Schr\"odinger equations in the form
\begin{align}
\left\{\begin{array}{ll}
 \left( i\pa_t+\frac{1}{2m_j}\pa_x^2 \right) u_j=N_j(u,\pa_x u), 
 & t>0,\ x \in \R,\ j=1,2,3,\\
 u_j(0,x)=\varphi_j(x), 
 & x\in \R,\ j=1,2,3
 \end{array}\right.
\label{eq}
\end{align}
where $u=(u_j)_{j=1,2,3}$ is a $\C^3$-valued unkonwn, 
$m_j$ is a positive constant and the nonlinear term $N_j$ satisfies
$$
 N_j(u,q)=O((|u|+|q|)^2) \quad \mbox{as} \quad (u,q) \to (0,0).
$$
We assume that $\varphi=(\varphi_j)_{j=1,2,3}$ belongs to the Sobolev space 
$H^s(\R)$ with $s \ge 1$, which is defined by 
$H^s(\R)= 
\{\psi\,;\, \pa_x^k \psi \in L^2(\R)\ \mbox{for all}\ k\le s \}$ 
equipped with the norm
$$
 \|\psi\|_{H^s} =\sum_{k \le s} \|\pa_x^k \psi\|_{L^2}.
$$

A typical nonlinear Schr\"odinger system appearing in various 
physical settings is 
\begin{align}
\left\{\begin{array}{l}
 \left( i\pa_t+\frac{1}{2m_1}\Delta \right) u_1= \overline{u_1} u_2, \\
 \left( i\pa_t+\frac{1}{2m_2}\Delta \right) u_2= {u_1}^2, 
\end{array}\right.
 \quad t>0,\ x \in \R^n
\label{2_waves}
\end{align}
(see e.g., \cite{CC}, \cite{CCO} for physical background). 
What is interesting in \eqref{2_waves} is that the ratio of the masses 
can affect the large-time behavior of the solutions. 
In the case of $n=2$, Hayashi--Li--Naumkin \cite{HLN} obtained a 
small data global existence result for \eqref{2_waves} under the 
relation $m_2=2m_1$. 
The non-existence of usual scattering state is also proved when $m_2=2m_1$. 
On the other hand, when $m_2\ne 2m_1$, it is shown in \cite{HLN2} that 
there is a usual scattering state under some restriction of the data.  
Higher dimensional case ($n \ge 3$) is considered by 
Hayashi--Li--Ozawa \cite{HLO} from the viewpoint of small data scattering.  
Remark that the relation $m_2=2m_1$ is often called the mass resonance 
relation, which was first discovered in the study of nonlinear Klein-Gordon 
systems (see \cite{DFX}, \cite{KOS}, \cite{KawaSuna}, \cite{Suna0}, 
\cite{Suna01}, \cite{Suna02}, \cite{Taf}, \cite{Y_Tsu}, etc.). 
More recently, large data case is discussed by Hayashi--Ozawa--Tanaka 
\cite{HOT}. In particular, their result includes finite time blow-up of the 
negative energy solutions for \eqref{2_waves} under mass resonance in the 
case of $4\le n \le 6$. 
However, their approach relies on the so-called virial identity which requires 
that the initial data of the blowing-up solutions must be suitably large 
(whence it should be distinguished from small data blow-up; 
see Section~5 below for more detail). Also it seems difficult to 
generalize blow-up results of this type to the case where the 
nolinearity involves the derivatives of the unkonwns. 
Concerning small data blow-up for NLS, very few results are 
known so far and many interesting problems are left unsolved (even in the case 
of single equations without derivatives in the nonlinear terms). 
We refer the readers to \cite{HNS}, \cite{IkeWak}, \cite{Oh}, 
\cite{OY}, \cite{Sasa}, \cite{Suna1}, \cite{Suna2} etc. for more information 
and the related topics.

The aim of this paper is to give examples of small data blow-up 
for (\ref{eq}). More precisely, we will show that there exist $m_j$, $N_j$ 
and $\varphi_j$ with $\|\varphi\|_{H^s}=\eps$ such that 
the corresponding solution blows up in finite time no matter how small 
$\eps>0$ is. We will also specify the order of the lifespan with respect to 
$\eps$. What we intend here is to illustrate, 
by using a simple model, how the interplay between the mass resonance 
and the nonlinear structure can affect global behavior of the solution. 
Although our examples below are somewhat artificial, they 
will help us to develop the understanding for possible mechanisms of 
singularity formation in more general nonlinear Schr\"odinger systems.

\section{Main result}
In what follows, we always assume that the nonlinearity in \eqref{eq} is 
in the form 
\begin{align} 
 N_1=0, \quad 
 N_2={u_1}^2, \quad 
 N_3=(\pa_x u_3)^2 + Q(u_1,u_2) \frac{\exp(2m_3 u_3)}{2m_3},
 \label{nonlin}
\end{align}
where $Q(u_1,u_2)$ is either 
${u_2}^2$, $u_1 u_2$, $\overline{u_1} u_2$ or $|u_2|u_2$. 
The main result is as follows. \\

\begin{thm} \label{thm1}
(1)\ \ 
Let $Q={u_2}^2$ and assume $m_1:m_2:m_3=1:2:4$. 
Then, for any $\eps \in (0,1]$ and $s \ge 1$, 
there exists $\varphi \in H^s(\R)$ with 
$\|\varphi\|_{H^s} = \eps$ such that the corresponding solution 
$u$ for \eqref{eq} satisfies 
\begin{align}
 \lim_{t \to T_{\eps} -0} \|u(t,\cdot)\|_{H^s}=\infty
 \label{blow-up}
\end{align}
with $T_{\eps} \in (\kappa \eps^{-4}, K \eps^{-4})$, 
where $\kappa$ and $K$ are positive constants not depending on $\eps$. 
\\
\indent(2)\ \ 
Let $Q=u_1 u_2$ and assume  $m_1:m_2:m_3=1:2:3$. 
Then, for any $\eps \in (0,1]$ and $s \ge 1$, 
there exists $\varphi \in H^s(\R)$ with 
$\|\varphi\|_{H^s} = \eps$ such that the corresponding solution  $u$ 
for \eqref{eq} satisfies \eqref{blow-up} 
with $T_{\eps} \in (\kappa \eps^{-6}, K \eps^{-6})$, 
where $\kappa$ and $K$ are positive constants not depending on $\eps$.
\\
\indent(3)\ \ 
Let  $Q=\overline{u_1}u_2$ and assume $m_1:m_2:m_3=1:2:1$. 
Then, for any $\eps \in (0,1]$ and $s \ge 1$, 
there exists $\varphi \in H^s(\R)$ with 
$\|\varphi\|_{H^s} = \eps$ such that the corresponding solution  $u$ 
for \eqref{eq} satisfies \eqref{blow-up} 
with $T_{\eps} \in (\kappa \eps^{-6}, K \eps^{-6})$, 
where $\kappa$ and $K$ are positive constants not depending on $\eps$.
\\
\indent(4)\ \ 
Let $Q=|u_2| u_2$ and assume  $m_1:m_2:m_3=1:2:2$. 
Then, for any $\eps \in (0,1]$, 
there exists $\varphi \in H^1(\R)$ with 
$\|\varphi\|_{H^1} = \eps$ such that the corresponding solution $u$ 
for \eqref{eq} satisfies 
\begin{align*}
 \lim_{t \to T_{\eps} -0} \|u(t,\cdot)\|_{H^1}=\infty
\end{align*}
with $T_{\eps} \in (\kappa \eps^{-4}, K \eps^{-4})$, 
where $\kappa$ and $K$ are positive constants not depending on $\eps$.\\
\end{thm}

\begin{rmk}
For general $\varphi \in H^s$ with $\|\varphi\|_{H^s}=\eps$, it is not 
difficult to show a lower bound for $T_\eps$ of the same order in $\eps$ 
(that is to say, we can show $T_{\eps}\ge \kappa \eps^{-4}$ 
in the case of (1), for instance) if $\eps$ is small enough. 
The novelty of the above theorem is the upper bound for $T_{\eps}$. 
In particular, this tells us that the order of the lifespan is actually 
influenced by the choice of $Q$ and the ratio of the masses. \\
\end{rmk}

\begin{rmk}
The relation between the choice of $Q$ and the ratio of the masses in 
Theorem~\ref{thm1} is characterized by the following condition:
\begin{align}
Q(e^{im_1\theta}z_1, e^{im_2\theta}z_2) = e^{im_3\theta}Q(z_1,z_2), 
\quad \theta \in \R,\ z_1,z_2 \in \C.
\label{gauge}
\end{align}
Our approach does not work without this condition. \\
\end{rmk}
\medskip

We close  this section by explaining our strategy of the proof. By setting
\begin{align}
 \sigma(t,x)=1-\exp(-2m_3 u_3(t,x)),
 \label{Cole_Hopf}
\end{align}
we can rewrite the system (\ref{eq}) with (\ref{nonlin}) as 
\begin{align*}
 \left\{\begin{array}{l}
 \left( i\pa_t+\frac{1}{2m_1}\pa_x^2 \right) u_1=0,\\
 \left( i\pa_t+\frac{1}{2m_2}\pa_x^2 \right) u_2={u_1}^2, \\
 \left( i\pa_t+\frac{1}{2m_3}\pa_x^2 \right) \sigma=Q(u_1,u_2).
\end{array}\right.
\end{align*}
This kind of transformation is first introduced by Hopf \cite{Hopf} 
and Cole~\cite{Cole} for the Burgers equation, and \eqref{Cole_Hopf} is 
used effectively by Ozawa \cite{Oz1}, \cite{Oz2} in the study of 
the quadratic NLS in the form 
$i\pa_t u+\frac{1}{2}\Delta u= (\nabla u)^2$ (see also p.38 of \cite{Stra}). 
Note that \eqref{Cole_Hopf} can be rewritten as 
$$ 
 u_3(t,x)=\frac{-1}{2m_3}\log(1-\sigma(t,x))
$$
if $|\sigma(t,x)|<1$, where the branch of the logarithm is chosen so that 
$\log 1=0$. 
Our main task in the proof of Theorem~\ref{thm1} is to choose $\varphi$ appropriately so that
\begin{align}
  \sigma(T_{\eps},x^*)=1 
 \label{sigma}
\end{align}
holds at some point $x^* \in \R$ (while $\|\sigma(t,\cdot) \|_{L^{\infty}}<1$ 
for $t<T_{\eps}$). The mass resonace condition 
(or, equivalently, \eqref{gauge}) 
will play a crucial role in the proof of this amplification. 
Once \eqref{sigma} is verified, we have
$$
 \|u_3(t,\cdot)\|_{H^s} 
 \ge 
 C|u_3(t,x^*)| 
 = 
 \frac{C}{2m_3} |\log(1-\sigma(t,x^*))|
 \to \infty
$$
as $t \to T_{\eps} -0$ 
(while $\|u_3(t,\cdot)\|_{H^s} <\infty$ for $t<T_{\eps}$). 
Similar idea can be found in the paper by Yagdjian \cite{Yag}, 
where semilinear wave equations with time-periodic coefficients are considered 
(see also \cite{KY}, \cite{OT}). 
Remark that the amplification in \cite{Yag} is due to parametric resonance and 
the proof is based on the Floquet theory.\\

\section{Preliminaries}
In this section, we collect several identities and estimates which are 
useful in the subsequent sections. 
In what follows, we denote several positive constants by the same letter $C$, 
which may vary from one line to another. 

First we put $\mathcal{L}_m=i\pa_t+\frac{1}{2m}\pa_x^2$ and 
$\mathcal{J}_m(t)=x+\frac{it}{m}\pa_x$ 
for $m>0$. Then we have $[\pa_x, \mathcal{J}_m(t)]=1$ and 
$[\mathcal{L}_m, \pa_x]=[\mathcal{L}_m, \mathcal{J}_m(t)]=0$, 
where $[\cdot, \cdot]$ denotes the commutator, i.e., 
$[\mathcal{P},\mathcal{Q}]=\mathcal{PQ}-\mathcal{QP}$ for 
linear operators $\mathcal{P}$ and $\mathcal{Q}$. 
Also we can easily check that 
\begin{align} 
 \mathcal{J}_{2m}(t)(\phi \psi)
 =\frac{1}{2}\bigl\{
    \bigl( \mathcal{J}_m(t) \phi \bigr) \psi 
    + 
    \phi \bigl( \mathcal{J}_m (t)\psi \bigr)
  \bigr\},
 \label{Leibniz}
\end{align}
\begin{align} 
 \mathcal{J}_{3m}(t)(\phi \psi)
 = \frac{1}{3} \bigl\{
   \bigl( \mathcal{J}_{m}(t) \phi \bigr) \psi 
    + 
    2 \phi \bigl( \mathcal{J}_{2m} (t)\psi \bigr)
  \bigr\}
\label{Leibniz2}
\end{align}
and
\begin{align} 
 \mathcal{J}_{m}(t)(\overline{\phi} \psi)
 =-\bigl(\overline{ \mathcal{J}_{m}(t) \phi} \bigr) \psi 
    + 
    2 \overline{\phi} \bigl( \mathcal{J}_{2m} (t) \psi)
 \label{Leibniz3}
\end{align}
for smooth functions $\phi$ and $\psi$. Next we put 
$\mathcal{A}_m(t)=\mathcal{F}_m \mathcal{U}_m(t)^{-1}$, 
where $\mathcal{F}_m$ and $\mathcal{U}_m(t)$ are defined by 
$$
 \bigl( \mathcal{F}_m \phi \bigr)(\xi)
 =
 \sqrt{\frac{m}{2\pi}}\int_{\R}e^{-imy\xi} \phi(y)dy
$$
and
$$
 \bigl( \mathcal{U}_m(t)\phi \bigr)(x)= 
 \sqrt{\frac{m}{2\pi it}}\int_{\R}e^{im(x-y)^2/(2t)} \phi(y)dy,
$$
respectively. Note that $w(t,x)=(\mathcal{U}_m(t)\phi)(x)$ solves 
$$
 \mathcal{L}_m w=0, \quad w(0,x)=\phi(x).
$$
We also remark that 
$\pa_t \mathcal{A}_m(t)\phi=-i\mathcal{A}_m(t) \mathcal{L}_m\phi$.

\begin{lem} \label{lem1}
For a smooth function $f(t,x)$, we have
$$
 \|f(t)\|_{L^{\infty}}
 \ge 
 t^{-1/2}\|\mathcal{A}_m(t)f(t) \|_{L^{\infty}} 
 -
 Ct^{-3/4} \rho_m[f](t),
$$
where
\begin{align}
 \rho_m[f](t)=\|f(t,\cdot)\|_{H^1} + \|\mathcal{J}_m(t) f(t,\cdot)\|_{L^2}.
 \label{def_rho}
\end{align}
\end{lem}
\medskip

\noindent{\bf Proof:}\ 
By the relation $\mathcal{J}_m(t)=\mathcal{U}_m(t)x\mathcal{U}_m(t)^{-1}$, 
we have
\begin{align*}
 \|\mathcal{A}_m (t)f\|_{H^1} 
 &=
 \|\mathcal{F}_m \mathcal{U}_m(t)^{-1}f\|_{H^1} \\
 &\le 
 C\|(1+|x|)\mathcal{U}_m(t)^{-1} f\|_{L^2}\\
 &\le 
 C\rho_m[f](t). 
\end{align*}
Next we observe that $\mathcal{U}_m(t)$ can be decomposed into 
$\mathcal{M}_m(t) \mathcal{D}(t) \mathcal{F}_m \mathcal{M}_m(t)$, where 
$$
 \bigl( \mathcal{M}_m(t) \phi \bigr)(x)=e^{imx^2/(2t)} \phi(x),
$$
$$
  \bigl( \mathcal{D}(t)\phi \bigr)(x)
  =
  \frac{1}{\sqrt{it}} \phi\left(\frac{x}{t}\right).
$$
We also set 
$\mathcal{W}_m(t)=\mathcal{F}_m M_m(t) \mathcal{F}_m^{-1}$. 
Then we see that 
\begin{align*} 
 f
 &= 
 \mathcal{U}_m(t) \mathcal{U}_m(t)^{-1}f\\
 &=
 \mathcal{M}_m(t) \mathcal{D}(t) \mathcal{F}_m \mathcal{M}_m(t) 
 \cdot \mathcal{F}_m^{-1} \mathcal{A}_m(t)f\\
 &=
 \mathcal{M}_m(t) \mathcal{D}(t) \mathcal{W}_m(t) \mathcal{A}_m(t)f.
\end{align*}
From the inequalities 
$\|f\|_{L^{\infty}} \le C \|f\|_{L^2}^{1/2}\|\pa_xf\|_{L^2}^{1/2}$ 
and 
$|e^{imx^2/(2t)}-1|\le C t^{-1/2}|x|$, 
it follows that 
\begin{align}
 \|(\mathcal{W}_m(t) -1)f\|_{L^{\infty}} 
 &\le 
 C\|(\mathcal{W}_m(t) -1)f\|_{L^{2}}^{1/2}
  \|\pa_x(\mathcal{W}_m(t) -1)f\|_{L^{2}}^{1/2}
 \nonumber\\
 &\le
 C (Ct^{-1/2} \|f\|_{H^{1}})^{1/2}(C\|f\|_{H^{1}})^{1/2} 
 \nonumber\\
 &=
 C t^{-1/4} \|f\|_{H^{1}}.
 \label{est_W}
\end{align}
Consequently we have 
\begin{align*}
 \|f-\mathcal{M}_m(t) \mathcal{D}(t) \mathcal{A}_m(t)f\|_{L^{\infty}}
 &=
 \|
   \mathcal{M}_m(t) \mathcal{D}(t) (\mathcal{W}_m(t) -1)\mathcal{A}_m(t)f
 \|_{L^{\infty}}\\
 &\le
 t^{-1/2} \|(\mathcal{W}_m(t) -1)\mathcal{A}_m(t)f\|_{L^{\infty}}\\
 &\le
 Ct^{-3/4} \|\mathcal{A}_m(t)f\|_{H^1}\\
 &\le 
 Ct^{-3/4} \rho_m[f](t),
\end{align*}
whence 
\begin{align*}
 \|f\|_{L^{\infty}}
 &\ge 
 \|\mathcal{M}_m(t) \mathcal{D}(t) \mathcal{A}_m(t)f \|_{L^{\infty}} 
 -
 \|f-\mathcal{M}_m(t) \mathcal{D}(t) \mathcal{A}_m(t)f\|_{L^{\infty}}\\ 
 &\ge 
 t^{-1/2}\|\mathcal{A}_m(t)f \|_{L^{\infty}} 
 -
 Ct^{-3/4}\rho_m[f](t)
\end{align*}
as required. 
\qed\\

\begin{lem} \label{lem2}
Let $f(t,x)$ and $g(t,x)$ be smooth functions satisfying 
$\mathcal{L}_{2m} g=f^2$. 
We have 
\begin{align} 
 \rho_{2m,s}[g](t) 
 \le 
 \rho_{2m,s}[g](0) 
 +  
 C \int_0^t  \rho_{m,s}[f](\tau)^2 \frac{d\tau}{\tau^{1/2}},
 \label{est_rho}
\end{align} 
where $\rho_{m,s}[\,\cdot\,]$ is defined by 
\begin{align*}
 \rho_{m,s}[f](t)
 =
 \|f(t,\cdot)\|_{H^s} + \|\mathcal{J}_m(t) f(t,\cdot)\|_{H^{s-1}}
\end{align*}
for $s \ge 1$. Also we have \begin{align} 
 \left\|
   \pa_t \bigl( \mathcal{A}_{2m}(t)g(t) \bigr)
   -
    e^{-i3\pi/4} t^{-1/2} \bigl( \mathcal{A}_{m}(t)f(t) \bigr)^2
 \right\|_{L^{\infty}} 
 \le 
 Ct^{-3/4} \rho_m[f](t)^2,
 \label{pointwise}
\end{align} 
where $\rho_m[f]=\rho_{m,1}[f]$, as defined in \eqref{def_rho}.
\end{lem}
\medskip

\noindent{\bf Proof:}\ 
First we note that 
$\mathcal{J}_m(t)=\frac{it}{m}\mathcal{M}_m(t)\pa_x \mathcal{M}_m(t)^{-1}$,
which implies
\begin{align} 
 \|f\|_{L^{\infty}}
 &=
 \|\mathcal{M}_m(t)^{-1}f\|_{L^{\infty}}
 \nonumber\\
 &\le
 C\|\mathcal{M}_m(t)^{-1}f\|_{L^2}^{1/2}
  \|\pa_x \mathcal{M}_m(t)^{-1}f\|_{L^2}^{1/2}
 \nonumber\\
 &\le
 C\|f\|_{L^2}^{1/2}  (t^{-1}\|\mathcal{J}_m(t)f\|_{L^2})^{1/2}
 \nonumber\\
 &\le 
 C t^{-1/2}\rho_m[f](t).
 \label{est_Sobolev}
\end{align} 
Since $[\mathcal{L}_{2m},  \pa_x^j \mathcal{J}_{2m}(t)^k]=0$, 
the standard energy method yields
\begin{align} 
 \frac{d}{dt}\| \pa_x^j \mathcal{J}_{2m}(t)^{k} g(t,\cdot)\|_{L^2} 
 \le 
 \|\pa_x^j \mathcal{J}_{2m}(t)^{k} (f^2)\|_{L^2}
 \label{est_energy}
\end{align}
for $k=0,1$ and $j \le s-k$. 
From \eqref{Leibniz},  \eqref{est_Sobolev} and \eqref{est_energy} 
it follows that  
\begin{align*} 
 \frac{d}{dt}\rho_{2m,s}[g](t) 
 &= 
 \sum_{k=0}^1 \sum_{j=0}^{s-k} \frac{d}{dt}
   \|\pa_x^j \mathcal{J}_{2m}(t)^{k} g(t,\cdot)\|_{L^2} \\
 &\le 
 \sum_{k=0}^1 \sum_{j=0}^{s-k} 
  \|\pa_x^j \mathcal{J}_{2m}(t)^{k} (f^2)\|_{L^2}\\
 &\le 
 C \left(
     \sum_{j'=0}^{s-1} \|\pa_x^{j'} f\|_{L^{\infty}} 
   \right)
   \left(
     \sum_{k=0}^1\sum_{j''=0}^{s-k} 
     \| \pa_x^{j''} \mathcal{J}_{m}(t)^{k} f\|_{L^2} 
   \right)\\
 &\le 
 C t^{-1/2} \rho_{m,s}[f](t)^2.
\end{align*}
By integrating with respect to $t$, we obtain the desired estimate 
\eqref{est_rho}. 
To prove \eqref{pointwise}, we put 
$\alpha(t,\xi)=\bigl(\mathcal{A}_{m}(t)f(t,\cdot) \bigr)(\xi)$, 
$\beta(t,\xi)=\bigl(\mathcal{A}_{2m}(t)g(t,\cdot) \bigr)(\xi)$ 
and
$$
 R(t,\xi)=i\pa_t \beta(t,\xi) - e^{-i\pi/4} t^{-1/2} \alpha(t,\xi)^2.
$$
Note that $\|\alpha(t,\cdot)\|_{H^1} \le C \rho_m[f](t)$ and  
\begin{align*}
 i\pa_t \beta
 =&
 \mathcal{A}_{2m}(t) \mathcal{L}_{2m} g\\
 =&
 \mathcal{F}_{2m} \mathcal{U}_{2m}(t)^{-1}(f^2)\\
 =&
 \mathcal{W}_{2m}(t)^{-1} \mathcal{D}(t)^{-1} \mathcal{M}_{2m}(t)^{-1} 
 \bigl( \mathcal{M}_{m}(t) \mathcal{D}(t) \mathcal{W}_{m}(t) \alpha \bigr)^2\\
 =&
 \mathcal{W}_{2m}(t)^{-1} \mathcal{D}(t)^{-1}  
 \bigl( \mathcal{D}(t) \mathcal{W}_{m}(t) \alpha \bigr)^2\\
 =&
  e^{-i\pi/4} t^{-1/2}
 \mathcal{W}_{2m}(t)^{-1} \bigl( \mathcal{W}_{m}(t) \alpha \bigr)^2.
\end{align*}
With the help of \eqref{est_W}, we have 
\begin{align*}
 \|R(t,\cdot)\|_{L^{\infty}}
 &= 
 t^{-1/2} 
 \bigl\|
  \mathcal{W}_{2m}(t)^{-1} \bigl( \mathcal{W}_{m}(t) \alpha \bigr)^2-\alpha^2 
 \bigr\|_{L^{\infty}}\\
 &\le 
 t^{-1/2} \Bigl(
 \bigl\|
   (\mathcal{W}_{2m}(t)^{-1}-1) \bigl( \mathcal{W}_{m}(t) \alpha \bigr)^2
 \bigr\|_{L^{\infty}}
 +
  \bigl\|
   \bigl( \mathcal{W}_{m}(t)\alpha \bigr)^2  - \alpha^2 \bigr)
  \bigr\|_{L^{\infty}}
 \Bigr)\\
 &\le 
 t^{-1/2} \Bigl(C t^{-1/4} 
  \bigl\| \bigl( \mathcal{W}_{m}(t) \alpha \bigr)^2 \bigr\|_{H^1}
 +
 \|(\mathcal{W}_{m}(t)+1)\alpha \|_{L^{\infty}}
 \|(\mathcal{W}_{m}(t)- 1) \alpha \|_{L^{\infty} }
 \Bigr) \\
 &\le 
 t^{-1/2} \Bigl(
 C t^{-1/4} 
 \| \mathcal{W}_{m}(t) \alpha \|_{L^{\infty}}
 \| \mathcal{W}_{m}(t) \alpha \|_{H^1}
 +
 C\|\alpha\|_{H^1} \cdot  C t^{-1/4} \| \alpha \|_{H^1}
 \Bigr) \\
&\le 
  C t^{-3/4} \| \alpha (t,\cdot)\|_{H^1}^2\\
&\le 
  C t^{-3/4} \rho_m[f](t)^2,
\end{align*}
which yields \eqref{pointwise}. 
\qed\\

\begin{rmk}
The above argument can be generalized as follows: 
Let $v_1$, $v_2$, $v_3$ be smooth fuctions of $(t,x)$ satisfying 
$\mathcal{L}_{m_3}v_3=Q(v_1,v_2)$, where $Q:\C \times \C \to \C$ 
satisfies \eqref{gauge} and 
$$
  Q(\lambda z_1,\lambda z_2)=\lambda^2 Q(z_1,z_2), 
 \quad \lambda >0,\ z_1,z_2 \in \C.
$$
Then we have 
$$
 i\pa_t \alpha_3=t^{-1/2} \tilde{Q}(\alpha_1, \alpha_2)+R,
$$
where 
\begin{align*}
 &\alpha_j(t,\xi)=\bigl(\mathcal{A}_{m_j}(t)v_j(t,\cdot) \bigr)(\xi), 
  \quad j=1,2,3,\\
 & \tilde{Q}(\alpha_1,\alpha_2)
   =
   e^{i\pi/4} Q(e^{-i\pi/4} \alpha_1, e^{-i\pi/4} \alpha_2),\\
 & R(t,\xi)
  =
   t^{-1/2} 
 \bigl\{
  \mathcal{W}_{m_3}(t)^{-1} 
  \tilde{Q}( \mathcal{W}_{m_1}(t) \alpha_1, \mathcal{W}_{m_2}(t) \alpha_2 )
  -
  \tilde{Q}( \alpha_1, \alpha_2 )
 \bigr\}.
\end{align*}
\end{rmk}

\section{Proof of Theorem~\ref{thm1}}
This section is devoted to the proof of Theorem~\ref{thm1}. 
Since the essential idea is the same, we consider the case of (1) in detail 
and only the outline of the proof will be given for the other cases. 
 
In what follows, we fix $s\ge 1$ and $\psi:\R\to \C$ which satisfies 
$\|\mathcal{F}_{m_1}^{-1}\psi\|_{H^s}=1$ and 
$\|(1+|x|)^{s-1}\pa_x \psi \|_{L^2}<\infty$. 
Let $v=(v_j(t,x))_{j=1,2,3}$ be the solution to 
\begin{align*}
 \left\{\begin{array}{l}
 \mathcal{L}_{m_1}  v_1=0,\\
 \mathcal{L}_{m_2} v_2={v_1}^2, \\
 \mathcal{L}_{m_3} v_3={v_2}^2,
\end{array}\right.
\quad t>0, \ \ x\in \R
\end{align*}
with the initial condition 
\begin{align}
 \left\{\begin{array}{l}
 v_1(0,x)=\eps (\mathcal{F}_{m_1}^{-1}\psi)(x),\\
 v_2(0,x)=0, \\
 v_3(0,x)=0.
 \end{array}\right.
 \label{initial_data}
\end{align}
We have the following. 

\begin{lem} \label{lem3}
Let $v$ be as above. 
Under the assumption $m_1:m_2:m_3=1:2:4$, 
there exist positive constants $\kappa$ and $K$, independent of $\eps$, 
such that 
\begin{align}
  \sup_{0\le t \le \kappa \eps^{-4}} \|v_3(t,\cdot)\|_{L^{\infty}}<1
 \label{est_upper}
\end{align}
and 
\begin{align}
  \|v_3(K \eps^{-4},\cdot)\|_{L^{\infty}}>1.
 \label{est_lower}
\end{align}
\end{lem}
\medskip

Before turning to the proof of Lemma~\ref{lem3}, we show that 
(1) of Theoroem~\ref{thm1} is derived from this lemma:  We set 
$T_{\eps} =\sup  \{T>0 \, ;\, \|v_3(t,\cdot)\|_{L^{\infty}}< 1
 \ \mbox{for}\ 0\le t<T \}$. 
Then \eqref{est_upper} and \eqref{est_lower} imply 
$\kappa \eps^{-4} < T_{\eps} < K \eps^{-4}$. 
Also, since the function $\R\ni x \mapsto |v_3(T_{\eps},x)|$ is continuous, 
we can choose $x^* \in \R$ such that 
$$
 |v_3(T_{\eps},x^*)|=\|v_3(T_{\eps},\cdot)\|_{L^{\infty}}=1.
$$
Now we take $\theta \in \R$ so that $v_3(T_{\eps},x^*)=e^{im_3\theta}$, 
and we set 
$$
 \varphi_1(x)=\eps e^{-im_1 \theta} (\mathcal{F}_{m_1}^{-1}\psi)(x), 
 \quad 
 \varphi_2(x)=\varphi_3(x)=0.
$$
Then, by the uniqueness of solutions to \eqref{eq}, we have
\begin{align*}
  &u_1(t,x)=e^{-im_1 \theta} v_1(t,x),\\
  &u_2(t,x)=e^{-im_2 \theta} v_2(t,x),\\
  &u_3(t,x)=\frac{-1}{2m_3} \log \bigl(1-e^{-im_3\theta} v_3(t,x) \bigr),
\end{align*}
which is a desired blowing-up solution. \qed \\

Now we are going to prove Lemma~\ref{lem3}. 
First we show \eqref{est_upper}. 
We put $\rho_{j,s}(t)=\rho_{m_j,s}[v_j](t)$ for $j=1,2,3$. 
By \eqref{est_rho}, we have 
\begin{align}
 \rho_{1,s}(t) \le \rho_{1,s}(0)=C\eps,
 \label{est_rho1}
\end{align}
\begin{align}
 \rho_{2,s}(t)
 \le 
 0 + 
 C\int_{0}^{t} \rho_{1,s}(\tau)^2  \frac{d\tau}{\tau^{1/2}} 
 \le
 C \eps^2 t^{1/2},
 \label{est_rho2}
\end{align}
and 
\begin{align}
 \rho_{3,s}(t)
 \le 
 0
 + 
 C\int_{0}^{t}  \rho_{2,s}(\tau)^2 \frac{d\tau}{\tau^{1/2}} 
 \le
 C \eps^4 t^{3/2}.
 \label{est_rho3}
\end{align}
From \eqref{est_Sobolev} and \eqref{est_rho3} it follows that  
$$
 \|v_3(t,\cdot) \|_{L^{\infty}} 
 \le  
 Ct^{-1/2} \rho_{3,s}(t) 
 \le 
 C\eps^4 t
 \le 
 C\kappa
$$
for $t \le \kappa \eps^{-4}$. 
By choosing $\kappa$ so small that $C\kappa<1$, 
we arrive at \eqref{est_upper}. 
Next we turn to the proof of \eqref{est_lower}. 
We put $\alpha_j(t,\xi)=(A_{m_j}(t) v_j(t,\cdot))(\xi)$ 
and 
$\rho_j(t)=\rho_{m_j,1}[v_j](t)$ 
for $j=1,2,3$. 
Since $\pa_t \alpha_1=-i\mathcal{A}_{m_1}(t)\mathcal{L}_{m_1}v_1=0$, 
we have 
\begin{align} 
 \alpha_1(t,\xi)
 =
 \alpha_1(0,\xi)
 =
 \eps \psi(\xi).
 \label{alpha1}
\end{align}
Also, it follows from \eqref{pointwise} that 
\begin{align*} 
 \left| 
  \pa_t{\alpha_{j+1}(t,\xi)} - 
  e^{-i3\pi/4} t^{-1/2} \alpha_j(t,\xi)^2 
 \right|
 \le 
  C t^{-3/4} \rho_j(t)^2
\end{align*}
for $j=1,2$. By \eqref{est_rho1} and \eqref{alpha1}, we have 
\begin{align} 
 |\alpha_2(t,\xi)- 2e^{-i3\pi/4}\eps^2 \psi(\xi)^2  t^{1/2}|
 &\le 
 |\alpha_2(1,\xi)- 2e^{-i3\pi/4} \eps^2 \psi(\xi)^2 | 
 +
 C \int_1^t \rho_1(\tau)^2  \frac{d\tau}{ \tau^{3/4}}
 \nonumber\\ 
 &\le 
 C\eps^2 + C\eps^2\int_1^t \frac{d\tau}{ \tau^{3/4}}
 \nonumber\\
 &\le 
 C\eps^2 t^{1/4}
 \label{asymp_alpha2}
\end{align}
for $t\ge 1$. As for $\alpha_3$, it follows from \eqref{est_rho2} that
\begin{align*} 
 \left|
 \alpha_3(t,\xi)- 
 e^{-i3\pi/4} \int_1^t (\alpha_2(\tau,\xi))^2 \frac{d\tau}{ \tau^{1/2}}
 \right|
 &\le 
 |\alpha_3(1,\xi)| + C\int_1^t \rho_2(\tau)^2  \frac{d\tau}{ \tau^{3/4}}\\
 &\le 
 C\eps^4 + C\eps^4\int_1^t \tau^{1/4} d\tau\\
 &\le 
 C\eps^4 t^{5/4}. 
\end{align*}
On the other hand, \eqref{asymp_alpha2} yields  
\begin{align*} 
 \left| 
  \int_1^t 
  (\alpha_2(\tau,\xi))^2 - (2e^{-i3\pi/4} \eps^2 \psi(\xi)^2  \tau^{1/2})^2
  \frac{d\tau}{\tau^{1/2}}
 \right|
 &\le
 \int_1^t
 C\eps^2 \tau^{1/4} \cdot C\eps^2 \tau^{1/2} \frac{d\tau}{\tau^{1/2}}
 \\
 &\le 
 C\eps^4 t^{5/4}.
\end{align*}
Summing up, we deduce that 
$$
 \Bigl|
 \alpha_3(t,\xi)  -\frac{8}{3}e^{-i9\pi/4} \eps^4 \psi(\xi)^4 t^{3/2} 
 \Bigr|
 \le 
 C\eps^4 t^{5/4}
$$
for $t\ge 1$. In particular, we obtain
\begin{align}
 \|\alpha_3(t,\cdot)\|_{L^{\infty}}
 &\ge 
  C^* \eps^4 t^{3/2}  - C\eps^4t^{5/4},
 \label{est_alpha3}
\end{align}
where $C^*= \frac{8}{3} \|\psi\|_{L^{\infty}}^4>0$. 
From \eqref{est_rho3}, \eqref{est_alpha3} and  Lemma~\ref{lem1} 
it follows that 
\begin{align*}
 \|v_3(t,\cdot)\|_{L^{\infty}}
 &\ge 
 t^{-1/2} \|\alpha_3(t,\cdot)\|_{L^{\infty}}  - C t^{-3/4} \rho_3(t) \\
 &\ge 
 C^* \eps^4 t  - C\eps^4t^{3/4}.
\end{align*}
By taking $K$ large enough, we  have 
\begin{align*}
 \|v_3(K\eps^{-4},\cdot)\|_{L^{\infty}} 
 &\ge 
 C^* K - C\eps K^{3/4} \\
 &\ge 
 (C^* K^{1/4} - C) K^{3/4} \\
 &>1,
\end{align*}
which completes the proof of \eqref{est_lower}. 
\qed\\

Finally, we give an outline of the proof of (2), (3), (4) of 
Theorem \ref{thm1}. In the case of (2), 
the problem is reduced to getting growth estimates for the solution 
$(v_j)_{j=1,2,3}$ to 
\begin{align*}
 \left\{\begin{array}{l}
 \mathcal{L}_{m_1}  v_1=0,\\
 \mathcal{L}_{m_2} v_2={v_1}^2, \\
 \mathcal{L}_{m_3} v_3={v_1}{v_2}
\end{array}\right.
\end{align*}
with the initial condition \eqref{initial_data}. 
Along the same line as the preceding argument, we can show that 
\begin{align} 
 \rho_{3,s}(t) \le C \eps^3 t
 \label{est_rho3_2}
\end{align}
and 
\begin{align*} 
 \| \alpha_3(t,\cdot)\|_{L^{\infty}} 
 \ge  
 2\eps^3 \|\psi\|_{L^{\infty}}^3 t - C\eps^3 t^{3/4}.
\end{align*}
Note that the identity \eqref{Leibniz2}, instead of \eqref{Leibniz}, 
plays the key role in the proof of \eqref{est_rho3_2}. 
By virtue of \eqref{est_Sobolev} and Lemma~\ref{lem1}, we have 
\begin{align*}
  \sup_{0\le t \le \kappa' \eps^{-6}} \|v_3(t,\cdot)\|_{L^{\infty}}<1
  \quad \mbox{and} \quad
  \|v_3(K' \eps^{-6},\cdot)\|_{L^{\infty}}>1
\end{align*}
with some positive constants $\kappa'$ and $K'$. It follows from these 
estimates that $T_{\eps} \in (\kappa'\eps^{-6}, K' \eps^{-6})$, 
which yields the desired conclusion. 
As for the proof of (3), we just have to replace 
\eqref{Leibniz2} with \eqref{Leibniz3} to obtain \eqref{est_rho3_2}. 
The proof of (4) is also similar: Just use
\begin{align}
 \|\mathcal{J}_{m_2}(t)(|v_2|v_2)\|_{L^2} 
 \le 
 C \|v_2\|_{L^{\infty}} \|\mathcal{J}_{m_2}(t)v_2 \|_{L^2} 
\label{Leibniz4}
\end{align}
in order to get the growth bound for $\rho_{3,1}(t)$. 
\qed\\

\section{Appendix: A quick review on blow-up of negative energy solutions
}

To make the difference between typical blow-up results and ours clearer, 
we will give a quick review on the proof of finite time blow-up for 
the 3-component NLS system 
\begin{align}
\left\{\begin{array}{l}
 \left( i\pa_t+\frac{1}{2m_1}\Delta \right) u_1= \overline{u_2} u_3, \\
 \left( i\pa_t+\frac{1}{2m_2}\Delta \right) u_2= \overline{u_1} u_3, \\
 \left( i\pa_t+\frac{1}{2m_3}\Delta \right) u_3= u_1 u_2, 
 \end{array}\right.
 \quad t>0,\ x \in \R^n,
\label{3_waves}
\end{align}
under the assumptions $E[u(0)]<0$, $m_3=m_1+m_2$ and $4 \le n \le 6$, 
where the energy $E[\,\cdot\,]$ is defined by 
$$
 E[\psi]=\sum_{j=1}^{3}\frac{1}{2m_j}\|\nabla \psi_j\|_{L^2}^2 
+2 \realpart \int_{\R^n} \psi_1(x) \psi_2(x) \overline{\psi_3(x)} dx 
$$
for $\psi=(\psi_j)_{j=1,2,3}$. 
Note that the 2-component system \eqref{2_waves} can be regarded 
as a degenerate case of \eqref{3_waves}, and the relation $m_3=m_1+m_2$ should 
be interpreted as the mass resonance relation associated with \eqref{3_waves}. 

The core of the proof is that the following three identities hold 
(cf. \cite{G}, \cite{M_Tsu}, etc.): 
\begin{align*}
&\frac{d}{dt} E[u(t)]=0, \\
&\frac{d}{dt} \sum_{j=1}^{3}m_j\|x u_j(t)\|_{L^2}^2 
 = 
 2V[u(t)] -  2(m_1+m_2-m_3)
 \imagpart \int_{\R^n} |x|^2 u_1(t,x) u_2(t,x) \overline{u_3(t,x)} dx,\\
&\frac{d}{dt} V[u(t)]=
 \frac{n}{2} E[u(t)] 
 +\frac{4-n}{2} \sum_{j=1}^{3}\frac{1}{2m_j}\|\nabla u_j(t)\|_{L^2}^2,
\end{align*}
where $V[\,\cdot\,]$ is defined by 
$$
 V[\psi] =\sum_{j=1}^{3}
 \imagpart \int_{\R^n} \overline{\psi_j(x)} x\cdot\nabla \psi_j(x)dx. 
$$
Once these identities are obtained, we can easily see that 
$$
 \sum_{j=1}^{3}m_j\|x u_j(t)\|_{L^2}^2 \le 
 \sum_{j=1}^{3}m_j\|x u_j(0)\|_{L^2}^2 + 2V[u(0)]t +\frac{n}{2}E[u(0)]t^2
 <0
$$
for sufficiently large $t$. 
This contradiction implies the non-existence of global solutions to 
\eqref{3_waves} in $H^1(\R^n) \cap L^2(\R^n; |x|^2dx)$ when 
$E[u(0)]<0$, $m_3=m_1+m_2$ and $n \ge 4$, 
while the local existence for \eqref{3_waves} can be shown 
when $n \le 6$, which comes from $p+1\le \frac{2n}{n-2}$ with $p=2$ 
(see \cite{HOT} for the detail).

We remark that $E[u(0)]<0$ implies $u(0)$ cannot be arbitrarily small, because
$$
 E[\eps \psi] =\eps^2 
 \left( 
 \sum_{j=1}^{3}\frac{1}{2m_j}\|\nabla \psi_j\|_{L^2}^2 
 + 2 \eps  \realpart \int_{\R^n} \psi_1(x) \psi_2(x) \overline{\psi_3(x)} dx 
 \right) 
 >0
$$
if $\psi \ne 0$ and $\eps>0$ is small enough. In fact, we can show 
the global existence of solutions to \eqref{3_waves} if the data are 
suitably small in $H^1(\R^n) \cap L^2(\R^n; |x|^2dx)$ when $n=4$  
(see \cite{HOT}). 
In this sense, our small data blow-up result presented in Theorem~\ref{thm1} 
should be distinguished from this kind  of ``large data" blow-up.

\medskip
\subsection*{Acknowledgments}
The first author (T.O.) is partially supported by Grant-in-Aid for Scientific 
Research~(A) (No.21244010), JSPS. 
The second author (H.S.) is partially supported by Grant-in-Aid for Young 
Scientists~(B) (No.~22740089), JSPS.


\end{document}